\documentclass[11pt]{amsart}
\usepackage{amsmath}
\usepackage{graphicx}
\usepackage[pdftex, pdfstartview=FitH]{hyperref}

\input xy
\xyoption{all}
\theoremstyle{definition}
\newtheorem{definition}{Definition}[section]

\newtheorem*{definition*}{Definition}
\newtheorem{question}{Question}

\theoremstyle{plain}
\newtheorem{theorem}[definition]{Theorem}
\newtheorem{corollary}[definition]{Corollary}
\newtheorem{lemma}[definition]{Lemma}
\newtheorem{proposition}[definition]{Proposition}

\newtheorem*{theorema*}{Theorem A}
\newtheorem*{theoremb*}{Theorem B}
\newtheorem*{theoremc*}{Theorem C}
\newtheorem*{theoremd*}{Theorem D}
\newtheorem*{lemma*}{Lemma}
\newtheorem*{proposition*}{Proposition}

\numberwithin{equation}{section}

\newcommand{\Proof}{\noindent {\it Proof.\ }}
\newcommand{\kernel}{{\rm Ker}\,}
\newcommand{\dotdot}{{\, . \, . \,}}
\newcommand{\hor}{{\mathcal H}}
\newcommand{\R}{{\mathbb R}}
\newcommand{\Hyp}{{\mathbb H}^2}

\newcommand{\Iso}{{\rm Iso}({\Hyp})}
\newcommand{\nmx}[1]{\mbox{$\|#1\|$}}
\newcommand{\grass}{{\widetilde{\rm G}}_2({\mathbb R}^4)}
\newcommand{\bx}{{\bf x}}
\newcommand{\bX}{{\bf X}}
\newcommand{\by}{{\bf y}}

\newcommand{\bu}{{\bf u}}
\newcommand{\bv}{{\bf v}}
\newcommand{\bw}{{\bf w}}

\newcommand{\bq}{{\bf q}}

\newcommand{\bi}{{\bf i}}
\newcommand{\bj}{{\bf j}}
\newcommand{\bk}{{\bf k}}
\newcommand{\bone}{{\bf 1}}
\newcommand{\bxi}{\mbox{\boldmath $\xi$}}

\newcommand{\haken}{\mathbin{\hbox to 6pt{%
                 \vrule height0.4pt width5pt depth0pt
                 \kern-.37pt
                 \vrule height6pt width0.4pt depth0pt\hss}}}

\title{Finsler surfaces with prescribed geodesics}

\author{J.C. \'Alvarez Paiva}
\address{J.C. \'Alvarez Paiva, U.M.R. CNRS 8524
         U.F.R. de Math\'ematiques, 59655 Villeneuve d'Ascq C\'edex, France}
\email{alvarez@math.univ-lille1.fr}

\author{G. Berck}
\address{G. Berck, Institut de Math\'ematiques, Universit\'e de Fribourg, Chemin du mus\'ee 23, CH-1700 Fribourg, Switzerland}
\email{gautier.berck@unifr.ch}
\thanks{Gautier Berck was partially funded by the Schweizerischer Nationalfonds grants SNF 200020-121506/1 and
 SNF PP002-114715/1 }
\keywords{Finsler manifold, geodesics, Crofton formula, path geometry}

\subjclass{53C65, 53C60, 53C22}

\begin{document}%
%%%%%%%%%%%%%%%%%

%%%%%%%%%%%%%%%%%%%%%%ABSTRACT%%%%%%%%%%%%%%%%%%%%%%%%%%%%%%%%%%%%%%%%%%
\begin{abstract}
We construct all Finsler metrics on the two-sphere for which geodesics are circles and show that any 
(reversible) path geometry on a two-dimensional manifold is locally the system of geodesics of a Finsler metric.
\end{abstract}
%%%%%%%%%%%%%%%%%%%%%%%%%%%%%%%%%%%%%%%%%%%%%%%%%%%%%%%%%%%%%%%%%%%%%%%%
\maketitle
%%%%%%%%%%%%%%%%%%%%%%%%%%%%%%%%%%%

%%%%%%%%%%%%%%%%%%%%%%
\section{Introduction}
%%%%%%%%%%%%%%%%%%%%%%
We consider the following variation of Hilbert's fourth problem:
{\sl to construct all Finsler metrics on the two-sphere whose geodesics are circles.}
Not great circles as in Hilbert's problem, just circles. Our primary goal is to
understand how to construct Finsler surfaces with prescribed geodesics. This question,
addressed by Busemann \cite{Busemann:1939, Busemann:1955, Busemann:1957},
Ambartzumian \cite{Ambartzumian:1976}, Alexander \cite{Alexander:1978},
Matsumoto \cite{Matsumoto:1989}, and Arcostanzo \cite{Arcostanzo:1994} will be completely
solved in Section~3. Our secondary goal is to construct all path geometries on the
two-sphere all of whose paths are circles.

We stress that we work with the classical definition of (reversible) Finsler metrics:

\begin{definition*}
A Finsler metric on a manifold $M$ is a continuous function $F : TM \rightarrow
[0 \dotdot \infty)$ that is homogeneous of degree one, smooth outside the zero section
and satisfies the {\sl quadratic convexity} condition: at every tangent space $T_xM$
the Hessian of $F^2(x,\cdot)$ (computed using any affine coordinates on $T_xM$) at
any nonzero tangent vector is a positive-definite quadratic form.
\end{definition*}

Indeed, in studying inverse problems in Finsler geometry our predecessors have either
weakened the condition of quadratic convexity to the condition that the restriction of $F$
to each tangent space be a norm (\cite{Arcostanzo:1994}), considered wider generalizations of 
Finsler metrics such as $G$-spaces (\cite{Busemann:1939, Busemann:1955, Busemann:1957}), or 
considered integrands that are homogeneous of degree one and not necessarily defined in all 
tangent directions (\cite{Matsumoto:1989}). If we insist on the quadratic convexity condition
it is because without it we cannot speak of the geodesic flow of the Finsler metric nor define important 
invariants such as the flag curvature (see \cite{Bao-Chern-Shen:2000}).

The paper consists of two parts each describing a geometric construction.

{\noindent \sl Part I. Finsler metrics associated to a family of hypersurfaces.}
Consider a manifold $M$ together with a family of hypersurfaces such that through
every point $x \in M$ and every tangent hyperplane $\zeta_x$ in $T_x M$ there is a unique
hypersurface passing through $x$ tangent to $\zeta_x$. Let us assume that this family is
parameterized by a smooth manifold $\Gamma$ and that the projection that sends a tangent
hyperplane $\zeta_x$ to the unique hypersurface that is tangent to it is a submersion from
the space $PT^* M$ of contact elements on $M$ to the parameter manifold $\Gamma$.

\begin{theorema*}\label{First-construction}
If $\mu$ is a smooth positive measure on $\Gamma$, there is a (unique) Finsler metric
$F : TM \rightarrow [0 \dotdot \infty)$ such that for any piecewise smooth curve $c$ on
$M$ we have the Crofton-type formula
\begin{equation}\label{Crofton-formula}
\int \!\! F(\dot{c}(t)) \ dt = \int_{\gamma \in \Gamma} \!\! \#(\gamma \cap c)
\ d\mu(\gamma) .
\end{equation}
\end{theorema*}

In general, the authors do not know what is the precise relationship between the family of
hypersurfaces and the Finsler metrics associated to it. However, in two dimensions the
relationship is simple enough:

\begin{theoremb*}\label{2D-case}
If $M$ is a two-dimensional manifold, the Finsler metrics defined by Eq.~(\ref{Crofton-formula}) are precisely
those whose geodesics are the curves of the family parameterized by $\Gamma$.
\end{theoremb*}

In particular, we will conclude in Section~3 that {\sl any path geometry in a two-dimensional manifold 
is locally the system of geodesics of a Finsler metric.} 

An interesting remark is that while the construction in Theorem~A extends to families of cooriented hypersurfaces 
such as the family of horospheres in hyperbolic space, Theorem~B does not. Indeed, if we allow different curves 
of $\Gamma$ to be tangent on the condition that their coorientations at the point of tangency be different,
then not only do those curves fail to be geodesics of the Finsler metric defined by Eq.~(\ref{Crofton-formula}), but
metrics associated to different measures may have different unparameterized geodesics. We shall demonstrate this in 
Section~4 by taking $\Gamma$ to be the family of horocycles in the hyperbolic plane and using Eq.~(\ref{Crofton-formula})
to construct two Riemannian metrics --- the hyperbolic metric and a metric conformal to it --- with different sets of 
unparameterized geodesics.

\noindent {\sl Part II. Circular path geometries on the two-sphere.}
These are families of circles on $S^2$ such that for every (tangent) line-element there
is a unique circle tangent to it. The standard example is the family of great circles on
the sphere.

The construction of all circular path geometries turns out to be quite simple and elegant.
Let us identify $\R^4$ with the algebra of quaternions and let $\pi : S^3 \rightarrow S^2$
be the Hopf fibration $\bq \mapsto \bq \bi \bar{\bq}$. Remark that any submanifold
$\Sigma$ of the Grassmannian $\grass$ of oriented two-planes in $\R^4$ engenders a family
of circles $\{\pi(L \cap S^3) : L \in \Sigma \}$ on the two-sphere.

Identifying, as usual, the Grassmannian $\grass$ with the product of two-spheres
$S_{-}^2 \times S_{+}^2$ of radius $1/\sqrt{2}$ in $\Lambda^2(\R^4) \cong \R^6$
(see Section~4 or Gluck and Warner~\cite{Gluck-Warner:1983} for details), the submanifolds
of $\grass$ that correspond to circular path geometries are very easy to construct:

\begin{theoremc*}
Let $\kappa : S_{+}^2 \rightarrow \R$ be a smooth function. The graph of the map
$f_\kappa : S_{+}^2 \rightarrow S_{-}^2$,
$$
f_\kappa (\bx) := \left(\frac{\kappa(\bx)}{\sqrt{2 + 2\kappa(\bx)^2}}, 0,
\frac{1}{\sqrt{2 + 2\kappa(\bx)^2}} \right),
$$
seen as a submanifold of the Grassmannian $\grass$, engenders a circular path geometry
if and only if $\kappa$ is odd and satisfies $\nmx{\nabla \kappa } < 1 + \kappa^2$. Moreover, every circular
path geometry can be obtained this way.
\end{theoremc*}

Let us denote by $\gamma_\kappa(\bx)$ ($\bx \in S_+^2$) the oriented circle on the two-sphere obtained by projecting down the intersection of
the three-sphere with the oriented two-plane $(f_\kappa(\bx),\bx) \in S_{-}^2 \times S_{+}^2$ by means of the Hopf fibration. Theorems~B and~C immediately imply the following simple description of all Finsler metrics on $S^2$ whose geodesics are circles.

\begin{theoremd*}
Let $\kappa$ be smooth odd function on $S_{+}^2$ such that $\nmx{\nabla \kappa} < 1 + \kappa^2$ and let $\mu$ be 
a smooth positive measure  on $S_{+}^2$. The Crofton-type formula
\begin{equation}
\int \!\! F(\dot{c}(t)) \ dt = \int_{\bx \in S_{+}^2} \!\! \#(\gamma_\kappa(\bx) \cap c) \ d\mu(\bx) 
\end{equation}
defines a Finsler metric on the sphere whose geodesics are circles. Moreover, every Finsler metric whose geodesics are 
circles can be obtained this way.
\end{theoremd*}

%%%%%%%%%%%%%%%%%%%%%%%%%%%%%%%%%%%%%%%%%%%%%%%%%%%%%%%%%%%%%%%%%
\section{Finsler metrics associated to a family of hypersurfaces}
%%%%%%%%%%%%%%%%%%%%%%%%%%%%%%%%%%%%%%%%%%%%%%%%%%%%%%%%%%%%%%%%%

In this section we show how to associate Finsler metrics to certain
families of hypersurfaces which are described by the class of double fibrations
\begin{equation}\label{Double-fibration}
\xymatrix{
      &  {S^* M} \ar[dl]_{\pi_1}\ar[dr]^{\pi_2} &   \\
      M  &     &  \Gamma ,}
\end{equation}
where $S^*M$ is the bundle of cooriented contact elements on $M$,
$\pi_1 : S^*M \rightarrow M$ is the canonical projection, $\Gamma$ is an oriented
manifold of the same dimension as $M$, $\pi_2 : S^* M \rightarrow \Gamma$ is a second
{\sl Legendrian fibration\/} (i.e., the fibers of $\pi_2$ are Legendrian submanifolds of
$S^*M$ provided with its canonical contact structure), and the map
$\pi_1 \times \pi_2 : S^* M \rightarrow M \times \Gamma$ is an embedding.

For each $\gamma \in \Gamma$, $\hat{\gamma} := \pi_1(\pi_2^{-1}\{\gamma\})$ is a
cooriented hypersurface. Notice that given a cooriented contact element $\xi \in S^*M$,
the hypersurface corresponding to $\pi_2(\xi) \in \Gamma$ is the unique hypersurface
from the family that is tangent to it and has a matching coorientation.

Notice that the description of families of non-cooriented hypersurfaces such as those appearing in 
the introduction can be obtained as a special case of~(\ref{Double-fibration}) by taking double covers. 

\noindent {\bf Example.}
The double fibration corresponding to all cooriented hyperplanes in $\R^n$ is
\begin{equation}\label{Standard-double-fibration}
\xymatrix{
      &  {S^{n-1} \times \R^n} \ar[dl]_{\pi_1}\ar[dr]^{\pi_2} &   \\
      \R^n  &     &  S^{n-1} \times \R ,}
\end{equation}
where the fibrations are given by $\pi_1(\xi,x) = x$ and
$\pi_2(\xi,x) = (\xi, \xi \cdot x)$.

In order to establish Theorem~A, and its generalization to the cooriented setting, we shall quickly review the basic
ideas around Gelfand transforms, Crofton formulas, and cosine transforms developed in~\cite{Alvarez-Fernandes:2007}.
In what follows {\sl we identify the bundle of cooriented contact elements on $M$, $S^* M$, with a hypersurface of 
$T^*M$ such that $S_x^* M$ is a star-shaped hypersurface in $T_x^* M$ for each $x \in M$.}

Given the double fibration~(\ref{Double-fibration}) and a volume form $\Omega$ on $\Gamma$, the Gelfand transform of
the volume density $|\Omega|$ is the function
$$
F := \pi_{1*}(\pi_2^* |\Omega|) : TM \longrightarrow \R
$$
whose value at the tangent vector $v \in T_xM$ is computed as follows: (1) construct a volume density on the fiber $\pi_1^{-1}\{x\} = S^*_x M$ by contracting $\pi_2^* |\Omega|$ at each point $\xi \in S^*_x M$ with any vector $\tilde{v} \in T_{\xi} S^* M$ such that
$D\pi_1(\tilde{v}) = v$; (2) integrate this volume density over $S^*_x M$.

\begin{lemma}\label{Gelfand-transform}
Given the double fibration~(\ref{Double-fibration}) and a volume form $\Omega$ on $\Gamma$, the Gelfand transform 
$F := \pi_{1*}(\pi_2^* |\Omega|)$ of the volume density $|\Omega|$ is homogeneous of degree one, smooth and positive
outside the zero section, and satisfies the Crofton-type formula
$$
\int \!\! F(\dot{c}(t)) \, dt = \int_{\gamma \in \Gamma} \!\! \#(\hat{\gamma} \cap c)
\, |\Omega|
$$
for any piecewise smooth curve $c$ on $M$.
\end{lemma}

\Proof The homogeneity, smoothness, and positivity all follow immediately from the construction. The Crofton-type formula is easily established
once we remark that $\pi_2^* |\Omega|$ defines a measure on the manifold $\pi_1^{-1}(c)$ whose dimension is the same as that of $\Gamma$ and hence, by the simplest case of the co-area formula,
\begin{eqnarray*}
\int_c \pi_{1*}(\pi_2^* |\Omega|) = \int_{\pi_1^{-1}(c)} \pi_2^* |\Omega| &=&
\int_{\gamma \in \Gamma} \! \#(\pi_2^{-1}\{\gamma\} \cap  \pi_1^{-1}(c)) \, |\Omega| \\
 &=&
\int_{\gamma \in \Gamma} \!\! \#(\hat{\gamma} \cap c) \, |\Omega| \, . \hfill \qed
\end{eqnarray*}

\noindent {\bf Example.}
Let us consider the case of $\R^n$ and its family of cooriented hyperplanes where the
double fibration is given in~(\ref{Standard-double-fibration}). If we write the volume
form  $\Omega$ on $S^{n-1} \times \R$ as $ \Omega = \omega \wedge dt$, where $\omega$ is
a volume form on $S^{n-1}$, then it is easy to see that
\begin{equation}\label{cosine-transform}
F(x,v) = \int_{\xi \in S^{n-1}} |\xi \cdot v|\, \omega .
\end{equation}

It is known since the work of Pogorelov on Hilbert's fourth problem
(see \cite{Pogorelov:1979}) that for any choice of volume form $\Omega$ the function $F$
in~(\ref{cosine-transform}) is a Finsler metric all of whose geodesics are straight lines.
For the more general double fibration~(\ref{Double-fibration}), we can say nothing about the
geodesics of $F$, but we can show that $F$ is a Finsler metric. 

\begin{theorem}\label{Crofton-construction}
If $\Omega$ is a volume form on $\Gamma$, then $F := \pi_{1*}(\pi_2^* |\Omega|)$ is a
Finsler metric on $M$ which, moreover, satisfies the Crofton-type formula
$$
\int \!\! F(\dot{c}(t)) \, dt = \int_{\gamma \in \Gamma} \!\! \#(\hat{\gamma} \cap c)
\, |\Omega|
$$
for any piecewise smooth curve $c$ on $M$.
\end{theorem}

In order to prove this theorem we first establish a formula for $F$ that
generalizes~(\ref{cosine-transform}). 

\begin{lemma}
For every point $x \in M$ there exists a volume form on $\omega_x$ defined on $S^*_x M$
such that the restriction of the function $F := \pi_{1*}(\pi_2^* |\Omega|)$ to the tangent
space $T_x M$ is given by the formula
\begin{equation}\label{Cosine-transform-analogue}
F(x,v) = \int_{\xi \in S^*_x M} |\xi(v)| \, \omega_x .
\end{equation}
\end{lemma}

\Proof
The key remark is that for every covector $\xi \in S^*_x M$ there exists a non-vanishing
$(n-1)$-form $\omega_{(x,\xi)}$ on $T_{(x,\xi)} S^*M$ such that
$$
\left(\pi_2^* \Omega \right)_{(x,\xi)} = \pi_1^* \xi \wedge \omega_{(x,\xi)} .
$$

Notice that $\kernel \xi$ is the tangent space at $\pi_1(\xi)$ of the hypersurface
corresponding to $\pi_2(\xi) \in \Gamma$. In other words,
$\kernel \xi = D\pi_1 (\kernel D_{(x,\xi)} \pi_2)$, or
$\kernel \pi_1^* \xi = \kernel D_{(x,\xi)}\pi_1 \oplus \kernel D_{(x,\xi)} \pi_2$.

It is now easy to verify that $\pi_1^* \xi \wedge \pi_2^* \Omega = 0$ and hence, by
Cartan's lemma, $\pi_1^* \xi$ is a linear factor of
$\left(\pi_2^* \Omega \right)_{(x,\xi)}$. Notice that while the form
$\omega_{(x,\xi)}$ is not completely determined, its restriction to
$\kernel D_{(x,\xi)}\pi_1 = T_\xi S^*_x M$ is. By performing this construction at every
point $\xi \in S^*_x M$ we obtain a non-vanishing $(n-1)$-form on $S^*_x M$ which we
shall call $\omega_x$. Clearly $\omega_x$ depends smoothly on the parameter $x$.

If $v \in T_x M$ is a tangent vector, the volume density constructed on $S^*_x M$
in order to push forward $\pi_2^* |\Omega|$ is just
$(x,\xi) \mapsto |\xi(v)||\omega_x|$. Since $\omega_x$ never vanishes, we may choose a
suitable orientation for $S^*_x M$ and conclude that
$$
F(x,v) = \pi_{1*}(\pi_2^* |\Omega|) (x,v) = \int_{S^*_x M} |\xi(v)|\omega_x .
$$
\qed

\noindent{\it Proof of Theorem~\ref{Crofton-construction}.\ }
In view of Lemma~\ref{Gelfand-transform} we just need to
verify the quadratic convexity condition. It is here that the relation between Gelfand transforms and Crofton
formulas proves useful through Eq.~(\ref{Cosine-transform-analogue}). Indeed, the quadratic convexity of norms 
defined by the cosine transform of smooth positive measures is a classical result which probably goes back to Blaschke, 
however we could not resist the temptation of presenting the following short proof.

Using the homogeneity of $F$ the problem reduces to showing that if $(v_1,\dots, v_n)$ are affine coordinates 
on $T_x M$, the Hessian of $F(x,\cdot)$ at any nonzero tangent vector $v$ is a positive semi-definite quadratic 
form whose kernel is the line spanned by $v$. Letting $(\xi_1, \dots, \xi_n)$ be the dual coordinates in $T^*_x M$, a 
smattering of distribution theory yields that

\begin{eqnarray*}
\frac{\partial}{\partial v_i} \int_{\xi \in S^*_x M} |\xi(v)|\, \omega_x
&=& \int_{\xi \in S^*_x M} \xi_i \, {\rm sign}(\xi(v)) \, \omega_x , \\
\frac{\partial^2}{\partial v_i \partial v_j} \int_{\xi \in S^*_x M} |\xi(v)|\, \omega_x
&=& 2 \int_{\xi \in S^*_x M} \xi_i \xi_j \, \delta(\xi(v)) \, \omega_x .\\
\end{eqnarray*}
From this formula, it follows that the quadratic form
$$
(D^2_v F)(w) = 2\int_{\xi \in S^*_x M} (\xi(w))^2 \, \delta(\xi(v)) \, \omega_x ,
$$
is positive semi-definite and its kernel is spanned by $v$. \qed

%%%%%%%%%%%%%%%%%%%%%%%%%%%%%%%%%%%%%%%%%%%%%%%%%%%%%%
\section{Metrizations of path geometries on surfaces }
%%%%%%%%%%%%%%%%%%%%%%%%%%%%%%%%%%%%%%%%%%%%%%%%%%%%%%

Intuitively, a path geometry on a surface $M$ is a family of curves such that for each
line-element in $PT^*M$ there is a unique curve tangent to it. As in the previous section,
it is more convenient to work ``upstairs" in the space of contact elements of $M$.

\begin{definition}
A {\sl path geometry\/} on a two-dimensional manifold $M$ is a smooth Legendrian
foliation of its space of contact elements $PT^*M$ whose leaves are transverse to
the canonical projection $\pi_1 : PT^*M \rightarrow M$. If $l$ is
a leaf of the path geometry, its projection $\pi_1(l)$ is called a {\sl path.}
\end{definition}

We consider the problem of determining whether a path geometry is the system of geodesics
of a Finsler metric. Since the Hopf-Rinow theorem (trivially) extends to Finsler metrics,
the following example of Howard~\cite{Howard:2005} shows, among other things, that there
are relatively simple path geometries which cannot be the geodesics of a Finsler metric.

\begin{theorem}[Howard~\cite{Howard:2005}]
There exists a complete, torsion-free affine connection on the two-torus $T^2$ such that
locally its geodesics can be made to correspond to lines on the plane, but, globally,
given any point $p \in T^2$, there is a second point $q$ such that there is no unbroken
geodesic joining $p$ and $q$.
\end{theorem}

Nevertheless we have the following result: 

\begin{theorem}\label{Metrizable-path-geometries}
If the space of leaves of a path geometry on a surface $M$ is a smooth manifold $\Gamma$
and the natural projection $\pi_2 : PT^*M \rightarrow \Gamma$ is a smooth fibration, then
for each smooth positive measure $\mu$ on $\Gamma$ there is a Finsler metric
$F : TM \rightarrow [0 \dotdot \infty)$ that
satisfies the Crofton-type formula
\begin{equation}\label{Crofton-for-surfaces}
\int \!\! F(\dot{c}(t)) \ dt = \int_{\gamma \in \Gamma} \!\! \#(\gamma \cap c)
\ d\mu(\gamma) ,
\end{equation}
and whose geodesics are the paths of the given path geometry. Moreover, any Finsler metric whose space
of geodesics is prescribed by $\Gamma$ can be constructed this way.
\end{theorem}

Locally every path geometry satisfies the hypotheses of
Theorem~\ref{Metrizable-path-geometries}. Indeed, by a result of
Whitehead~\cite{Whitehead:1932} around any point in a manifold $M$ provided with a
path geometry, there exists a {\sl convex neighborhood\/} $U$ where any two distinct
points in $U$ are joined by unique path segment lying entirely in $U$. This allows
us to identify the space of paths lying in $U$ with the (smooth) quotient of the
manifold
$$
\{(x,y) \in \partial U \times \partial U : x \neq y \}
$$
by the equivalence relation $(x,y) \cong (y,x)$. As a consequence,
Theorem~\ref{Metrizable-path-geometries} has the following corollary:

\begin{corollary}
Any path geometry on a surface $M$ is locally the system of geodesics of a Finsler metric.
\end{corollary}

\noindent{\it Proof of Theorem~\ref{Metrizable-path-geometries}.\ }
The Finsler metric $F := \pi_{1*}(\pi_2^* \, \mu)$, constructed as in
Theorem~\ref{Crofton-construction}, automatically satisfies the Crofton-type
formula~(\ref{Crofton-for-surfaces}). In turn, the following classical argument shows
that the existence of a Crofton-type formula implies that paths are geodesics. 

Let $p$ and $q$ be points in a convex neighborhood $U$ and let $\gamma_{pq}$ be unique path segment 
that lies entirely in $U$ and joins the two points. Notice  that if $\gamma' \in \Gamma$
is a path different from the one containing $\gamma_{pq}$, then any connected component of $\gamma' \cap U$ intersects 
$\gamma_{pq}$ at most once. Indeed, if it intersects it at two points we would have two distinct path segments inside $U$ joining
the same two points, which contradicts the convexity of $U$. Furthermore, if $\sigma$ is a piecewise smooth curve in $U$ 
joining $p$ and $q$, any connected component of $\gamma' \cap U$ must intersect $\sigma$ at least as many times as it
intersects $\gamma_{pq}$. Using the Crofton-type formula, we have
$$
{\rm length}(\sigma) = \int_{\gamma' \in \Gamma} \#(\sigma \cap \gamma') \, d\mu
\geq
  \int_{\gamma' \in \Gamma} \#(\gamma_{pq} \cap \gamma') \, d\mu
= {\rm length}(\gamma_{pq})
$$
and, therefore, path segments are locally length-minimizing.

Let us now prove the converse. Assume that the family of geodesics of a Finsler metric $F$ is precisely our family of paths and
let us find a measure $\mu$ on $\Gamma$ for which the Crofton-type formula~(\ref{Crofton-for-surfaces}) holds. 

Let $\tilde{\Gamma}$ be the double cover of $\Gamma$ which parameterizes the family of paths with co-orientation. The natural 
projection $\pi_2 : PT^*M \rightarrow \Gamma$ lifts to a projection $\tilde{\pi}_2 : S^*M \rightarrow \Gamma$, where $S^*M$ is the
bundle of unit covectors. If $\omega$ is the standard symplectic form on $T^*M$ and $i : S^*M \rightarrow T^*M$ is the inclusion,
there is a unique area-form $\tilde{\omega}$ on $\tilde{\Gamma}$ such that the equality 
$i^* \omega = \tilde{\pi}_2^* \, \tilde{\omega}$ holds on $S^*M$. By the general Crofton formula in Finsler spaces 
(Theorem~5.2 in \cite{Alvarez-Berck:2006}), for any smooth curve $c$
$$
\int \!\! F(\dot{c}(t)) \ dt = \frac{1}{4} \int_{\gamma \in \tilde{\Gamma}} \!\! \#(\gamma \cap c) \, |\tilde{\omega}| .
$$
It follows that $(1/2)|\tilde{\omega}|$ descends to a smooth, positive measure $\mu$ on $\Gamma$ for which the Crofton-type 
formula~(\ref{Crofton-for-surfaces}) holds. 
\qed

%%%%%%%%%%%%%%%%%%%%%%%%%%%%%%%%%%%%%%%%%%%%%%%%%%%%%%%%%%%%%%%%%
\section{Horocycles and Crofton formulas on the hyperbolic plane}
%%%%%%%%%%%%%%%%%%%%%%%%%%%%%%%%%%%%%%%%%%%%%%%%%%%%%%%%%%%%%%%%%

One remarkable feature of Theorem~\ref{Metrizable-path-geometries} is that different measures in 
Eq.~(\ref{Crofton-for-surfaces}) yield different Finsler metrics with the same geodesics. We shall now see in an
example that this is no longer the case for non-reversible path geometries where different cooriented 
paths can be tangent on the condition that their coorientations at the point of tangency be different. 

Consider $\R^3$ together with the Lorentzian form $[\bx,\by] = -x_1 y_1 - x_2 y_2 + x_3 y_3$. As is well-known,
the hyperbolic plane can be defined as the hypersurface
$$
\Hyp = \{\bx \in \R^3 : [\bx,\bx] = 1, x_3 > 0 \}
$$
with the induced Riemannian metric $g_\bx(\bv,\bw) = -[\bv,\bw]$. The unit tangent bundle of $\Hyp$ can be
identified with the submanifold
$$
S\Hyp = \{(\bx,\bv) : \bx \in \Hyp, [\bx,\bv] = 0, [\bv,\bv] = -1  \} \subset \R^3 \times \R^3 \ 
$$
and the family of horocycles is parameterized by the upper light-cone
$$
\hor = \{\bxi \in \R^3 : [\bxi,\bxi] = 0, \xi_3 > 0 \}:
$$
to each $\bxi \in \hor$ we associate the horocycle  $\{\bx \in \Hyp : [\bxi, \bx] = 1 \}$. 

The incidence relation associated to the family of horocycles is described by the double fibration
\begin{equation}\label{horocyclic-double-fibration}
\xymatrix{
      &  {S\Hyp} \ar[dl]_{\pi_1}\ar[dr]^{\pi_2} &   \\
      \Hyp  &     &  \hor ,}
\end{equation}
where the fibrations are given by $\pi_1(\bx,\bv) = \bx$ and
$\pi_2(\bx,\bv) = \bx + \bv$. Indeed, it is easy to verify that the horocycle associated to $\bx + \bv$ is precisely 
the horocycle that passes through $\bx$ and whose tangent is orthogonal to (and cooriented by) the unit vector $\bv$.

The construction of a measure on $\hor$ that is invariant under the induced action of isometry group $\Iso$ is also 
quite simple: the two-form $\left((d\xi_1 \wedge d\xi_2 \wedge d\xi_3) \haken \bxi \right)/[\bxi,\bxi]$ in $\R^3$ is manifestly 
invariant under Lorentz transformations and restricts to a well-defined volume form $\Omega$ on the upper light cone. Another
way to describe this volume form is as $(d\xi_1 \wedge d\xi_2)/\xi_3$.  

In the two propositions that follow we shall prove that the Finsler metrics associated to the measures $|\Omega|$ and $\xi_3 |\Omega| = |d\xi_1 \wedge d\xi_2|$ 
do not have the same geodesics.

\begin{proposition}
Consider the measures $\Omega$ and $\xi_3 |\Omega| = |d\xi_1 \wedge d\xi_2|$ on the space of horocycles $\hor$. The Finsler
metric $F = \pi_{1*}(\pi_2^* |\Omega|)$ is a positive multiple of the hyperbolic metric and the Finsler metric 
$\pi_{1*}(\pi_2^* \, \xi_3|\Omega|)$ is conformal to $F$ with conformal factor $x_3$. 
\end{proposition}

\Proof
The $\Iso$-equivariance of the double fibration~(\ref{horocyclic-double-fibration}) and the invariance of the measure $|\Omega|$, immediately
imply that $F$ is invariant and hence a positive multiple of the hyperbolic metric. On the other hand
$$
\pi_{1*}(\pi_2^* \, \xi_3|\Omega|) = \pi_{1*}((x_3 + v_3) \, \pi_2^*|\Omega|) = \pi_{1*}(x_3 \, \pi_2^* |\Omega|) + \pi_{1*}(v_3 \, \pi_2^* |\Omega|).
$$ 
Since $v_3$ is an odd function on $\pi_1^{-1}(x) = S_\bx \Hyp$, the second summand is zero  and 
$\pi_{1*}(\pi_2^* \, \xi_3|\Omega|) = x_3 F$.
\qed

\begin{proposition}
Two Finsler metrics that are conformal by a non-constant factor do not have the same geodesics.
\end{proposition}

\Proof
Let $F$ be a Finsler metric on a manifold $M$ and let $\nu$ be a non-constant, positive, smooth function on $M$. It is enough to
show that $F$ and $\nu F$ do not have the same geodesics on some open chart where $\nu$ is not constant and, therefore, we may work
with local coordinates. The Euler-Lagrange equations for $\nu F$ are
$$
\frac{d}{dt}\left(\nu \frac{\partial F}{\partial v}\right) - \frac{\partial (\nu F)}{\partial x} = 0 .
$$ 
Therefore, a path $x(t)$ that is a geodesic for both $F$ and $\nu F$ must satisfy the equations
$$
\frac{d}{dt}\left(\frac{\partial F}{\partial v}\right) - \frac{\partial F}{\partial x} = 0 \mbox{  and  }
\frac{d\nu}{dt}\frac{\partial F}{\partial v} - F \frac{\partial \nu}{\partial x} = 0.
$$
However, if we choose initial conditions $(x,v)$ such that $(\partial \nu / \partial x) (x) \neq 0$ and 
$(\partial \nu / \partial x) (x) \cdot v = 0$, we see that the second equation cannot be satisfied
unless $v = 0$. \qed

%%%%%%%%%%%%%%%%%%%%%%%%%%%%%%%%%%%%%%%%%%%%%%%%%%%%
\section{Circular path geometries on the two-sphere}
%%%%%%%%%%%%%%%%%%%%%%%%%%%%%%%%%%%%%%%%%%%%%%%%%%%%

Our aim in this section is to construct all path geometries in the two-sphere where the
paths are circles. In fact, such {\sl circular path geometries \/} will be identified
with a certain class of great circle fibrations of the three-sphere. The results of
Gluck and Warner~\cite{Gluck-Warner:1983} will then yield a simple construction.

We begin by using the Hopf fibration $\pi : S^3 \rightarrow S^2$,
$\bq \mapsto \bq \bi \bar{\bq}$, to establish a simple correspondence between oriented
circles on the two-sphere of purely imaginary quaternions of unit length and a class of
oriented great circles on the three-sphere of unit quaternions.

Let $(\bq, \kappa) \in S^3 \times \R$ and consider the oriented great circle
$C(\bq,\kappa) \subset S^3$ that passes through $\bq$ in the direction of
$\bq \bu$, where
$$
\bu = \frac{\kappa \bi}{\sqrt{1 + \kappa^2}} + \frac{\bk}{\sqrt{1 + \kappa^2}} .
$$

It is easy to see that the oriented great circles $C(-\bq,\kappa)$ and $C(\bq,\kappa)$
coincide and that $C(\bq,\kappa)\bi$ is the same great circle, but with
reversed orientation, as $C(\bq \bi, -\kappa)$.

\begin{proposition}\label{circle-construction}
The image $c(\bq, \kappa)$ of the great circle $C(\bq,\kappa)$ under the Hopf
fibration is the oriented circle on the two-sphere that passes through
$\pi(\bq) = \bq \bi \bar{\bq}$ in the direction of $\bq \bj \bar{\bq}$ and has geodesic
curvature $\kappa$. Moreover, for any oriented circle $c$ on the two-sphere there is a
unique great circle $C(\bq,\kappa)$ such that $c = c(\bq, \kappa)$.
\end{proposition}

\Proof
Simply put, we will identify circles in $S^2$ with their Frenet frames in $SO(3)$ and then
use the standard identification of the double cover of the group of rotations with the
group of unit quaternions.

If $\bx : (a \dotdot b) \rightarrow S^2$ is a regular curve and
$\bv(t) = \dot{\bx}(t)/\nmx{\bx(t)}$ is its normalized velocity vector, then the
matrix $\bX(t)$ whose columns are $\bx(t)$, $\bv(t)$, and $\bx(t) \times \bv(t)$
describes a curve in $SO(3)$. The Frenet equations take the form
$$
\dot{\bX}(t) = \bX(t)
\left( \begin{array}{ccc}
0 & -s(t) & 0 \\
s(t) & 0 & -s(t)\kappa(t) \\
0 & s(t)\kappa(t) & 0 \end{array} \right) ,
$$
where $s(t) = \nmx{\dot{\bx}(t)}$ is the speed of the curve and $\kappa(t)$ is its
geodesic curvature. Notice that the curve $\bx(t)$ describes a circle with constant speed
if and only if the matrix on the right is constant and, therefore, if and only if the
Frenet frame $X(t)$ is the left-translation of a one-parameter subgroup on $SO(3)$.

Consider now the map $\sigma : S^3 \rightarrow SO(3)$ that takes the unit quaternion
$\bq$ to the rotation $\sigma(\bq)(\bx) = \bq \bx \bar{\bq}$, $\bx \in S^2$. In matrix
form, $\sigma(\bq)$ is the matrix whose columns are the vectors $\bq \bi \bar{\bq}$,
$\bq \bj \bar{\bq}$, and $\bq \bk \bar{\bq}$. The map $\sigma$ is a
two-to-one covering ($\sigma(\bq) = \sigma(-\bq)$) and a group homomorphism. It
follows that $\sigma$ defines a bijection between great circles in $S^3$ and
left-translates of one-parameter subgroups in $SO(3)$. Moreover, an easy computation
shows that the image of a great circle $C(\bq,\kappa)$ under $\sigma$ satisfies the
equation
$$
\dot{\bX}(t) = \bX(t)
\left( \begin{array}{ccc}
0 & -2/\sqrt{1 + \kappa^2} & 0 \\
2/\sqrt{1 + \kappa^2} & 0 & -2\kappa(t)/\sqrt{1 + \kappa^2} \\
0 & 2\kappa(t)/\sqrt{1 + \kappa^2} & 0 \end{array} \right) .
$$
In other words, $\sigma(C(\bq,k))$ is the Frenet lift of a circle of curvature $\kappa$
described with constant speed $2/\sqrt{1 + \kappa^2}$. \qed

Given a circular path geometry on the two-sphere we may construct a foliation of $S^3$ by
oriented great circles as follows: given a point $\bq \in S^3$, consider the unique
(circular) path $c$ on the two-sphere that passes through $\bq \bi \bar{\bq} = \bx$ and is
tangent to the line spanned by $\bq \bj \bar{\bq} = \bv$ at this point. Orient $c$ so that
$\bv$ is its normalized tangent vector at $\bx$. If $\kappa$ is the geodesic curvature
of the oriented circle $c$, take $C(\bq,\kappa)$ to be the leaf of the the foliation that
passes through $\bq$. Another way of understanding this foliation is to remark that the
Frenet lifts of all circles in the path geometry, where each circle is oriented in two
different ways, form a foliation of $SO(3)$ which can be transferred to $S^3$ via the
double cover $\sigma : S^3 \rightarrow SO(3)$.

The (smooth)  great circle foliations of the three-sphere associated to paths geometries on the
two-sphere are easy to characterize. First, if we provide $S^3$ with the contact structure
such that the contact hyperplane at the point $\bq$ is the plane spanned by $\bq \bi$ and
$\bq \bk$, then all great circles in the foliation are Legendrian. Second, if $C$ is an
oriented great circle in the foliation and $\dot{\bq}$ is its velocity vector
at a point $\bq \in C$, then $\langle \dot{\bq},\bq\bk \rangle > 0$. Third, since we lift
every circle in the path geometry with its two possible orientations, if $C$ is an
oriented great circle in the foliation, then $C \bi$ {\sl with its orientation reversed \/}
is also in the foliation.

Nevertheless, let us forget for an instant the special character of these great circle
foliations and recall Gluck and Warner's simple characterization of (general) great circle
foliations of the three-sphere. The idea is to characterize all submanifolds of
the Grassmannian of oriented two-planes in $\R^4$ that induce great circle foliations
of $S^3$.

Using coordinates associated to the basis of $\Lambda^2 \R^4$ formed by the bivectors
\begin{eqnarray*}
&\bv_1 = (\bone \wedge \bi - \bj \wedge \bk)/2, \
\bv_2 = (\bone \wedge \bj + \bi \wedge \bk)/2, \
\bv_3 = (\bone \wedge \bk - \bi \wedge \bj)/2, \\
&\bv_4 = (\bone \wedge \bi + \bj \wedge \bk)/2, \
\bv_5 = (\bone \wedge \bj - \bi \wedge \bk)/2, \
\bv_6 = (\bone \wedge \bk + \bi \wedge \bj)/2,
\end{eqnarray*}
the set of simple bivectors of unit length is the product of the spheres
\begin{eqnarray*}
S_{-}^2 &=& \{(x_1,x_2,x_3,0,0,0) : x_1^2 + x_2^2 + x_3^2 = 1/2 \} \quad {\rm and} \\
S_{+}^2 &=& \{(0,0,0,x_4,x_5,x_6) : x_4^2 + x_5^2 + x_6^2 = 1/2 \} .
\end{eqnarray*}

Using this identification of $\grass$ with $S_{-}^2 \times S_{+}^2$, Gluck and Warner
proved the following

\begin{theorem}[\cite{Gluck-Warner:1983}]\label{Gluck-Warner-thm}
A submanifold of $\grass$ corresponds to a smooth fibration of the three-sphere
by oriented great circles if and only if it is the graph of a smooth map
$f : S_{\pm}^2 \rightarrow S_{\mp}^2$ with $\|df\| < 1$.
\end{theorem}

We are now ready to prove the main result in this section.

\begin{theoremc*}
Let $\kappa : S_{+}^2 \rightarrow \R$ be an smooth function. The graph of the map
$f : S_{+}^2 \rightarrow S_{-}^2$,
$$
f(\bx) := \left(\frac{\kappa(\bx)}{\sqrt{2 + 2\kappa(\bx)^2}}, 0,
\frac{1}{\sqrt{2 + 2\kappa(\bx)^2}} \right),
$$
seen as a submanifold of the Grassmannian $\grass$, engenders a circular path geometry
if and only if $\kappa$ is odd and satisfies $\nmx{\nabla \kappa} < 1 + \kappa^2$. Moreover, every circular
path geometry can be obtained this way.
\end{theoremc*}

\Proof
We just need to characterize the type of map $f$ appearing in
Theorem~\ref{Gluck-Warner-thm} that corresponds to the class of oriented great
circle fibrations of $S^3$ associated to circular path geometries on $S^2$.

In the fibrations we are considering, an oriented great circle $C$ passing through
a point $\bq$ does so in the direction of $\bq \bu$, where $\bu$ is of the form
$a\bi + b \bk$ with $a^2 + b^2 = 1$ and $b > 0$. Since the coordinates of the unit
bivector $(1/\sqrt{2}) \bq \wedge \bq \bu$ are easily computed to be of the form
$(1/\sqrt{2})(a,0,b,x,y,z)$, where $x^2 + y^2 + z^2 = 1$ we have that the map
corresponding to the fibration must a map $f : S_{+}^2 \rightarrow S_{-}^2$ of the form
$$
f(x) = (\alpha(\bx),0,\beta(\bx)), \quad \alpha(\bx)^2 + \beta(\bx)^2 = \frac{1}{2},
\quad \beta(\bx) > 0.
$$
Letting $\kappa(\bx) = \alpha(\bx)/\beta(\bx)$, we have that
$$
f(\bx) := \left(\frac{\kappa(\bx)}{\sqrt{2 + 2\kappa(\bx)^2}}, 0,
\frac{1}{\sqrt{2 + 2\kappa(\bx)^2}} \right).
$$
and the condition that $\|df\| < 1$ translates to $\nmx{\nabla \kappa } < 1 + \kappa^2$.

It remains to show that the function $\kappa$ is odd. To see this recall the symmetry
condition: if $C$ is an oriented great circle in the foliation associated to a path geometry,
then $C \bi$ with its orientation reversed is also in the foliation. In coordinates
$(x_1,\dots,x_6)$, the action of right multiplication by $\bi$ on the space of bivectors
is given by
$$
(x_1,x_2,x_3,x_4,x_5,x_6) \mapsto (x_1,-x_2,-x_3,x_4,x_5,x_6) .
$$
Therefore, we must have that $(-\alpha(\bx),0,\beta(\bx),-\bx)$ is in the graph of $f$
and so $\alpha$ is odd, $\beta$ is even, and $\kappa = \alpha/\beta$ is odd.
\qed

%%%%%%%%%%%%%%%%%%%%%%%
\section{Final Remarks}
%%%%%%%%%%%%%%%%%%%%%%%

It is natural to ask how many of the metrics constructed in Theorem~D are actually Riemannian. The answer is {\it not many.} Indeed,
it follows at once from a result of A.~G. Khovanskii~(\cite{Khovanskii:1980}) that if the geodesics of an affine connection on $S^2$
are circles, then there exists a diffeomorphism that sends the family of geodesics to the family of great circles. It follows
then from Beltrami's theorem that the only Riemannian metrics that can be constructed through Theorem~D are of constant curvature. The difference between Riemannian and Finsler geometry responsible for this phenomenon is that the exponential map of a Finsler metric is not $C^2$ on the zero section. 

In a similar vein, we may ask for the Riemannian version of Theorem~B: {\sl when is a two-dimensional path geometry locally the family of (unparameterized) geodesics of a Riemannian metric on $M$?} The preliminary problem of determining whether the paths are locally the (unparameterized) geodesics of an affine connection is classical and has been studied by Tresse, Cartan, Bol, and many others. We refer the reader to the book of Arnold~(\cite{Arnold:1988}) and the paper of Bryant~(\cite{Bryant:1997}) for two excellent modern expositions. Even when the family of paths is locally the family of geodesics of an affine connection, there are obstructions and it is only recently that substantial progress has been made by Bryant, Dunajski, and Eastwood~(\cite{Bryant-Dunajski-Eastwood:2008}) towards the solution of this problem. 

Another interesting problem is to ascertain up to what point Theorem~B extends to the non-reversible setting. We distinguish two different
questions.

\begin{question}
Given a two-dimensional reversible path geometry, construct all non-reversible Finsler metrics for which the paths
are geodesics.
\end{question}

\begin{question}
Is every non-reversible two-dimensional path geometry locally the family of geodesics of a non-reversible Finsler metric?
\end{question}

The challenge these questions pose stems from the inability of integral-geometric techniques such as cosine transforms and Crofton-type 
formulas to yield non-reversible metrics. 

The simplest example of a non-reversible path geometry on the plane is the family of circles of a fixed radius. This case has been elegantly solved by Tabachnikov in~\cite{Tabachnikov:2004}. In addition to the non-reversibility of this path geometry, one encounters an
additional phenomenon: the Finsler metrics whose geodesics are circles of a fixed radius can only be defined locally. Indeed, a version of the Hopf-Rinow theorem holds for non-reversible Finsler metrics (see Theorem~6.6.1 in \cite{Bao-Chern-Shen:2000}) and is a clear obstruction
to the existence of globally-defined Finsler metrics whose geodesics are circles of a fixed radius. Thus it seems that in this and other 
non-reversible problems it is better to work with the following generalization of Finsler metrics:

\begin{definition*}
A {\sl magnetic Lagrangian} on a manifold $M$ is a continuous function $L : TM \rightarrow \R$ that is positively
homogeneous of degree one, smooth outside the zero section and satisfies the following {\sl quadratic convexity condition:} 
the restriction of $L$ to each tangent space $T_xM$ is the support function of a quadratically-convex body in $T_x^*M$. 
\end{definition*}

The quintessential magnetic Lagrangian on the plane is 
$$
L(x,y;u,v) = \sqrt{u^2 + v^2} + \frac{k}{2}(xv-yu) ,
$$  
which models the motion of a particle under a constant magnetic field and whose extremals are circles of a fixed radius.

To conclude, we clarify the relationship between Theorem~D and Theorem~4.2 in~\cite{Alvarez-Duran:2002} where
\'Alvarez Paiva and Dur\'an construct a large class of examples of Finsler metrics on the sphere whose geodesics are circles. 

The idea in~\cite{Alvarez-Duran:2002} is quite simple: let us take a smooth positive measure $\mu$ on the space $S^{3*}$ of cooriented
totally geodesic two-spheres in $S^3$ and require that it be invariant under the circle action induced from the Hopf action on $S^3$. By
Pogorelov's solution of Hilbert's fourth problem, the Crofton-type formula   
$$
\int \! F(\dot{c}(t) \, dt = \int_{\xi \in S^{3*}} \! \#(\xi \cap c) \, d\mu(\xi)
$$
defines a Finsler metric that is invariant under the Hopf action and for which great circles are geodesics. If $\pi : S^3 \rightarrow S^2$
is the Hopf fibration, we construct a submersive Finsler metric on the two-sphere by requiring that the unit ball at $T_{\pi(\bq)}S^2$ be
the image under $D\pi$ of the unit ball at $T_\bq S^3$. We refer to \cite{Alvarez-Duran:2002} for the proof that the geodesics of this metric are circles. 

We wish to point out that this construction cannot give rise to every Finsler metric whose geodesics are circles. Indeed, the invariance of 
the measure $\mu$ under the induced Hopf action implies that this construction depends on a function of two variables, while
Theorem~D shows that the general construction depends on two functions of two variables.

\bibliography{../../paperbib}
\bibliographystyle{amsplain}
%%%%%%%%%%%%%%%
\end{document}